\documentclass[preprint]{elsarticle}
\usepackage{placeins}
\usepackage{amsmath}
\usepackage{amssymb}
\usepackage{tikz}
\usepackage{graphicx}
\usepackage{epstopdf}
\usepackage{natbib}

\begin{document}

\title{Stochastic resetting in the Kramers problem: a Monte Carlo approach \tnoteref{t1}}
\tnotetext[t1]{This work has been supported by the Spanish State Research Agency (AEI) and the European Regional Development Fund (ERDF, EU) under Project No.~PID2019-105554GB-I00.}

\author[1]{Julia Cantis\'{a}n \corref{cor1}}

\cortext[cor1]{Corresponding author}
\ead{julia.cantisan@urjc.es}

\author[1]{Jes\'{u}s M. Seoane}
\author[1,2]{Miguel A.F. Sanju\'{a}n}

\address[1]{Nonlinear Dynamics, Chaos and Complex Systems Group, Departamento de F\'{i}sica, Universidad Rey Juan Carlos \\ Tulip\'{a}n s/n, 28933 M\'{o}stoles, Madrid, Spain}

\address[2]{Department of Applied Informatics, Kaunas University of Technology \\ Studentu 50-415, Kaunas LT-51368, Lithuania}

 \begin{abstract}
 The theory of stochastic resetting asserts that restarting a search process at certain times may accelerate the finding of a target. In the case of a classical diffusing particle trapped in a potential well, stochastic resetting may decrease the escape times due to thermal fluctuations. Here, we numerically explore the Kramers problem for a cubic potential, which is the simplest potential with a escape. Both deterministic and Poisson resetting times are analyzed. We use a Monte Carlo approach, which is necessary for generic complex potentials, and we show that the optimal rate is related to the escape times distribution in the case without resetting. Furthermore, we find rates for which resetting is beneficial even if the resetting position is located on the contrary side of the escape.
 \end{abstract}


 \maketitle

\section{Introduction}

\begin{figure}[h]
	\begin{center}
		\includegraphics[width=0.45\textwidth]{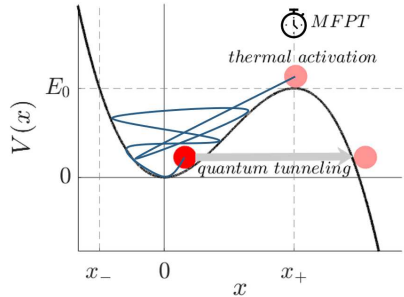}
	\end{center}
	\caption{Two mechanisms allow the particle to overcome the potential barrier: thermal activation and quantum tunneling. For a classical particle, after some time the noise induces escapes from the well. The time necessary for the particle to escape is related to the time needed to reach the top of the potential barrier, that is, the mean first passage time (MFPT).}
	\label{V_shape}
\end{figure}

In principle, a particle trapped in a potential well is bounded in that region of space forever. The bottom of the well corresponds to the stable state and under any small perturbation the system returns to that point. This is because the system has no energy to overcome the potential barrier. However, two mechanisms can change the particle's fate, allowing it to escape: thermal activation and quantum tunneling (see Fig. \ref{V_shape}). Both effects may be present and one may dominate over the other depending on the temperature of the system \cite{Grabert1984, Munakata1987, Sharifi1988, Vlasov2016}.

If we deal with a classical problem where quantum tunneling is not possible, it is the presence of fluctuations and dissipation (that arise naturally in the presence of a thermal bath, for instance) which eventually overcomes the potential barrier. In other words, there is a nonzero probability to escape from the potential well and if we wait for a sufficiently long time it is reasonable to expect that the particle can escape from it. The problem that determines the rate at which a Brownian particle escapes from a potential well is called the Kramers problem \cite{Kramers1940}. For a review of the problem one may read \cite{Melnikov1991} and for a recent review on the subject and an extension to field theory, see \cite{Berera2019}.

The time necessary for the particle to escape is related to the time needed to reach the top of the potential barrier. This time is usually referred to as mean first passage time (MFPT). We use the term \textit{mean} as we are dealing with a stochastic process and different trajectories escape at different times.

In this paper we mainly investigate a mechanism that affects the MFPT: stochastic resetting (SR). It is based on the idea that resetting to a given state (usually the initial state) may accelerate a search process, in our case, a escape process where the target is the top of the energy barrier. SR was first studied by Evans and Majumdar \cite{Evans2011}. They showed that the MFPT of a diffusing particle searching for a target in an infinite domain becomes finite when SR is implemented. Also, they showed that MFPT presents a minimum for a certain value of the resetting rate. The intuitive idea behind this is that by restarting the process we avoid the realizations that get too far from the target. Evans and Majumdar also showed that resetting creates a nonequilibrium steady state, in contrast to a normal diffusing particle for which the width of the Gaussian position distribution grows with time as $ \sqrt{t} $.

Some natural processes such as animal foraging \cite{Viswanathan2011} or human visual search \cite{Wolfe2004} utilize this mechanism to reduce the search time. Also, it is known that restarting computer algorithms may enhance their performance \cite{Lorenz2018}, as in simulated annealing \cite{Li2016}, which is used to locate the global minimum of a potential landscape. At low temperatures, the algorithm may get stuck in a local minimum for long periods of time and resetting gives the algorithm a second chance. This method has also been applied to the process of RNA polimerization \cite{Roldan2016} or, recently, to the Ising model \cite{Magoni2020}. For more applications, see \cite{Evans2020}.

In recent years, SR has been a hot topic for the research community. Several papers have dealt with multiple targets \cite{Bressloff2020}, multiple searchers \cite{Bhat2016}, resetting in non-diffusive processes such as Levy flights \cite{Kusmierz2014} or telegraphic processes \cite{Masoliver2019}. Different time-resetting distributions have been explored and also diffusion in some simple potentials (logarithmic \cite{Ray2020}; linear, quadratic, box potentials \cite{Ahmad2019}; piece-wise linear potential \cite{Singh2020}). Recently, SR was also experimentally studied using holographic optical tweezers for the resetting \cite{Besga2020, Tal-Friedman2020}.

The vast majority of the literature on SR is analytical, thus it is limited to more or less simple systems that are mathematically tractable. Here, we explore another approach using Monte Carlo simulations that can be used for complex systems that arise non closed-form expressions. We apply it to the problem of diffusion in a nonmonotonic potential that has one escape. We take into account the statistical fluctuations of this method, but we show that it is on agreement with some general analytic predictions. This numerical approach allows to perform a complete study of the problem. We explore not only the MFPT, but the FPT distribution. For example, we show that the optimal resetting rate, $ r^{*} $, is related to the FPT distribution for the system without resetting.

Finally, the article is organized as follows. In section \ref{Section_2} we present the problem without resetting. We motivate our numerical approach and show how the noise and the initial state affect the process. In sections \ref{Section_3} and \ref{Section_4} resetting at deterministic and Poisson times are analyzed. Finally, we provide concluding remarks at the end.

\section{Diffusion in the Kramers problem} \label{Section_2}

In a real world situation, diffusion can not always be considered to take place in an infinite domain. On the contrary, diffusion is usually subjected to a bias. These situations can be modeled as diffusion in a potential landscape, that is, the Kramers problem. The Kramers problem is studied in diverse areas such as atomic diffusion, where the potential accounts for the interatomic forces in a crystal lattice \cite{Vega2002}, or in chemical reactions, where the potential barrier is the activation energy needed to the break chemical bonds \cite{Dzhioev2011}. Another application is that of Josephson junctions, for which the simplest model consists of a classical particle in a tilted washboard potential in the presence of friction and thermal noise \cite{Costantini1999}.

We consider the diffusion of a classical particle trapped in the cubic potential associated to the Helmholtz oscillator \cite{Almendral2004}:
\begin{equation}
	V(x)= \frac{\alpha x^{2}}{2}- \frac{\beta x^{3}}{3},
	\label{V}
\end{equation}
where we take $\alpha=6$ and $ \beta=1 $, to obtain a potential shape as in Fig. \ref{V_shape}. It is a simple asymmetrical nonlinear potential with an escape on the right-hand side of the well and an infinite wall on the other side. It is used to study a diverse range of nonlinear phenomena such as \cite{Casado1994}. The bottom of the well is located at $ x=0 $ and the height of the potential barrier is $ E_{0}=\alpha^{3}/(6 \beta^{2})=36 $. The noiseless basin of attraction of the origin is delimited by $ V(x)=36 $, that is, by $ x_{+} $, the top of the potential barrier, and $ x_{-} $, the left-hand side position at the same energy level. In our case, we have $ x_{+}=\alpha / \beta=6 $ and $ x_{-}=-\alpha / (2 \beta)=-3 $.

Including the stochastic noise in Newton's second law, we get the Langevin equation for our process:
\begin{equation}
	\ddot{x}+ \eta \dot{x}+V'(x)= \xi(t),
	\label{Langevin_general}
\end{equation}
where $ \eta \dot{x}$ is the damping term and $ \xi(t) $ is a Gaussian white noise, so that $ \left\langle \xi(t) \right\rangle =0 $ and $ \left\langle \xi(t) \xi(t') \right\rangle = \varepsilon \delta(t-t') $, where $ \varepsilon $ accounts for the noise intensity. The diffusion constant, $ D $, is related to $ \varepsilon $ through Einstein's relation: $ \varepsilon= 2D=2 \eta k_{B}T $, where $ k_{B} $ is the Boltzmann constant and $ T $ is the temperature.

We take $ \eta=0.1 $ so that we remain in the range of light damping as $ \eta^{2} < \omega_{0}^{2} $, where $ \omega_{0} $ refers to the natural frequency of the noiseless system. Including the derivative of the potential we get the final equation:
\begin{equation}
	\ddot{x}+0.1 \dot{x}+6x-x^{2}= \xi(t).
	\label{Langevin}
\end{equation}

Throughout the whole paper we consider small noise intensities compared to the height of the potential barrier $ E_{0} $: $ \varepsilon \ll E_{0} $, to model realistic situations and to make sure that the trajectories do not escape immediately. For instance, for $ \varepsilon=1.8 $ the noise is a $ 5 \% $ of the height of the potential barrier: $\varepsilon/(\alpha^{3}/(6 \beta^{2}))=1.8/36=0.05.$

The damping and the noise terms are responsible for the system's instability and allow escapes from the potential well. On one hand the damping dissipates kinetic energy and slows down the particle hindering the return to the equilibrium point. On the other hand, the noise acts as a random force that pushes the particle away from the initial position. This combined effect results in a nonzero probability to escape from the potential well.

In contrast with other stochastic resetting problems that did not considered any escapes, our potential has one escape. In order to compute the MFPT, the position $ x_{+}=6 $ can be thought as a trap that absorbs the diffusing particle once it is hit for the first time. Under this criterion the MFPT is the same as the mean escape time. For higher noise and damping levels, we would have to consider the chance of the particle coming back once it has passed $ x_{+}=6 $. Simulations were not stopped until further positions were achieved and no come back was observed. Thus, for our parameters we are safe equating the MFPT and the mean escape time.

\subsection{Motivation: analytic versus Monte Carlo approach}

Traditionally, stochastic resetting problems involve the analytical computation of the steady state probability density and the MFPT. Some papers have included numerical results to show the good agreement between both approaches, but a complete numerical analysis of the problem is lacking. Besides, an analytic approach is not always possible for complex potentials and numerical approaches are needed to explore more realistic situations.

Here, we deal with a SR problem with a potential whose analytical solution is difficult to find because the calculations lead to a non closed-form expression as we show in the following paragraphs. Details of the numerical calculations and its limitations are also presented in this subsection.

Let us start with the analytical approach. We consider a particle whose equation of motion is described by Eq. \ref{Langevin}. This particle is reset to its initial state $ (x_{0}, \dot{x}_{0}) $ at different times drawn from a distribution with a constant rate, $ r $. For the first passage problem we would like to solve the related backward Chapman-Kolmogorov equation as we are interested in the probability of the particle to survive until time $ t $, starting from any initial position $x $, that is $ Q(x,t) $. Denoting $ Q_{0} \equiv Q(x_{0},t)$, we have:
\begin{equation}
	\frac{\partial Q}{\partial t} = D  \frac{\partial^{2} Q}{\partial x^{2}} -V'(x)  \frac{\partial Q}{\partial x} -rQ+rQ_{0},
	\label{ChapmanK}
\end{equation}
where the last two terms account for a loss and a gain of probability due to the resetting. The case with resetting can be derived from the case without resetting once it is calculated \cite{Reuveni2016}. Thus, for practical purposes, we solve firstly the problem without resetting:
\begin{equation}
	\frac{\partial Q}{\partial t} = D \frac{\partial^{2} Q}{\partial x^{2}} - V'(x)  \frac{\partial Q}{\partial x},
	\label{ChapmanK_V}
\end{equation}
with initial condition $ Q(x,0)=1 $ and boundary condition $ Q(6,t)=0 $. In order to solve this equation we take the Laplace transform:
\begin{equation}
	q(x,s)=\int_{0}^{\infty} e^{-st} Q(x,t) \,dt,
	\label{LaplaceT}
\end{equation}
which simplifies the problem to:
\begin{equation}
	D\frac{d^2 q}{dx^2}-V'(x)\frac{d q}{dx}-sq+1=0,
	\label{Heun}
\end{equation}
with $ q(6,s)=0 $. Notice that $ V'(x)=6x-x^{2} $ in our case. Hence, Eq. \ref{Heun} is a special case of the triconfluent Heun equation. The solution to this equation is given in terms of the Heun functions:
\begin{equation}
	\begin{split}
		q(x)&=A \times \mathrm{HeunT}(-10^{\frac{1}{3}} s,3,-3 \,10^{\frac{2}{3}},-\frac{10^{\frac{1}{3}} (x-3)}{3}) \\
		&+ B \times {\mathrm e}^{-\frac{10 x^{2} (x-9)}{27}} \times \mathrm{HeunT}(-10^{\frac{1}{3}} s,-3,-3 \,10^{\frac{2}{3}},\frac{10^{\frac{1}{3}} (x-3)}{3})+\frac{1}{s},
	\end{split}
\end{equation}
which are a type of advanced special functions \cite{Slavyanov2000}. There is no general closed-form expression for these type of functions as the coefficients of their Frobenius expansions obey recurrence relations of at least three terms. The theory of Heun functions is currently under study and has brought attention in the context of quantum control and engineering \cite{Ishkhanyan2018}.

For simpler potentials (logarithmic \cite{Ray2020}; linear, quadratic, box potentials \cite{Ahmad2019} or a piece-wise linear potential \cite{Singh2020}), Eq. \ref{Heun} is exactly solvable  as the coefficient of the spatial derivative is simpler.

For complex potentials, we may use a Monte Carlo (MC) approach to avoid the analytical difficulties. The drawback of this method is that it implies statistical fluctuations and it is subject to error. However, we show that the error can be kept small and that the results are on agreement with some general analytical predictions.

The idea behind the Monte Carlo method is that if we want to estimate the theoretical first passage time distribution (which is described by a distribution with  mean $ \mu $ and standard deviation $ \sigma $), we can take samples ($ X_{1}, X_{2},...X_{N} $), i.e., realizations of the stochastic process and average the passage times. By the law of large numbers we know that if there are sufficient samples the average must be close to the true value. Also, the Central Limit Theorem establishes that the mean of a sample with size $ N $ is a random variable with mean value  $ \tilde{\mu}=\mu $ and standard deviation $ \tilde{\sigma}=\sigma/ \sqrt{N} $ \cite{Reiter2008}. The standard deviation usually is not known, but it may be estimated as
\begin{equation}
	\sigma=\sqrt{\frac{1}{N} \sum^{N}_{i=1}(X_{i}-\tilde{\mu})^2}.
\end{equation}

Thus, the Monte Carlo error in our estimation of the MFPT is $ \sigma/ \sqrt{N} $. This is logical as the error decreases for very similar values and for a large number of realizations. However, the error goes with the square root of $ N $ and this means that increasing $ N $ to get one more decimal digit of accuracy requires a factor of $ 100 $ increase in the number of trials, and thus, more computation time and memory storage. Finally, the Central Limit Theorem may be used to establish confidence intervals. For a confidence interval of $ 95 \% $, we express our results of the MFPT as:
\begin{equation}
	MFPT=\tilde{\mu} \pm \frac{2 \sigma}{ \sqrt{N}}.
	\label{MC_error}
\end{equation}

This uncertainty due to the statistical fluctuations is expressed in the form of error bars in our graphs.

We numerically integrate Eq. \ref{Langevin} for $ 10^{4} $ trajectories, so that the error is kept small and computation times are manageable. We used the adaptive SOSRA algorithm in Julia \cite{Bezanson2017}, which is specific for problems with additive noise \cite{Rackauckas2020} and as options we chose reltol=abstol=$10^{-3}  $. We measured the time for which each realization reached $ x=6 $ and counted that as the first passage time. Simulations were not stopped until $ x>100 $ to ensure no comebacks happened. Different initial locations are considered but the initial velocity, $ \dot{x_{0}} $, is always considered to be zero.

\begin{figure}[h]
	\begin{center}
		\includegraphics[width=0.45\textwidth ]{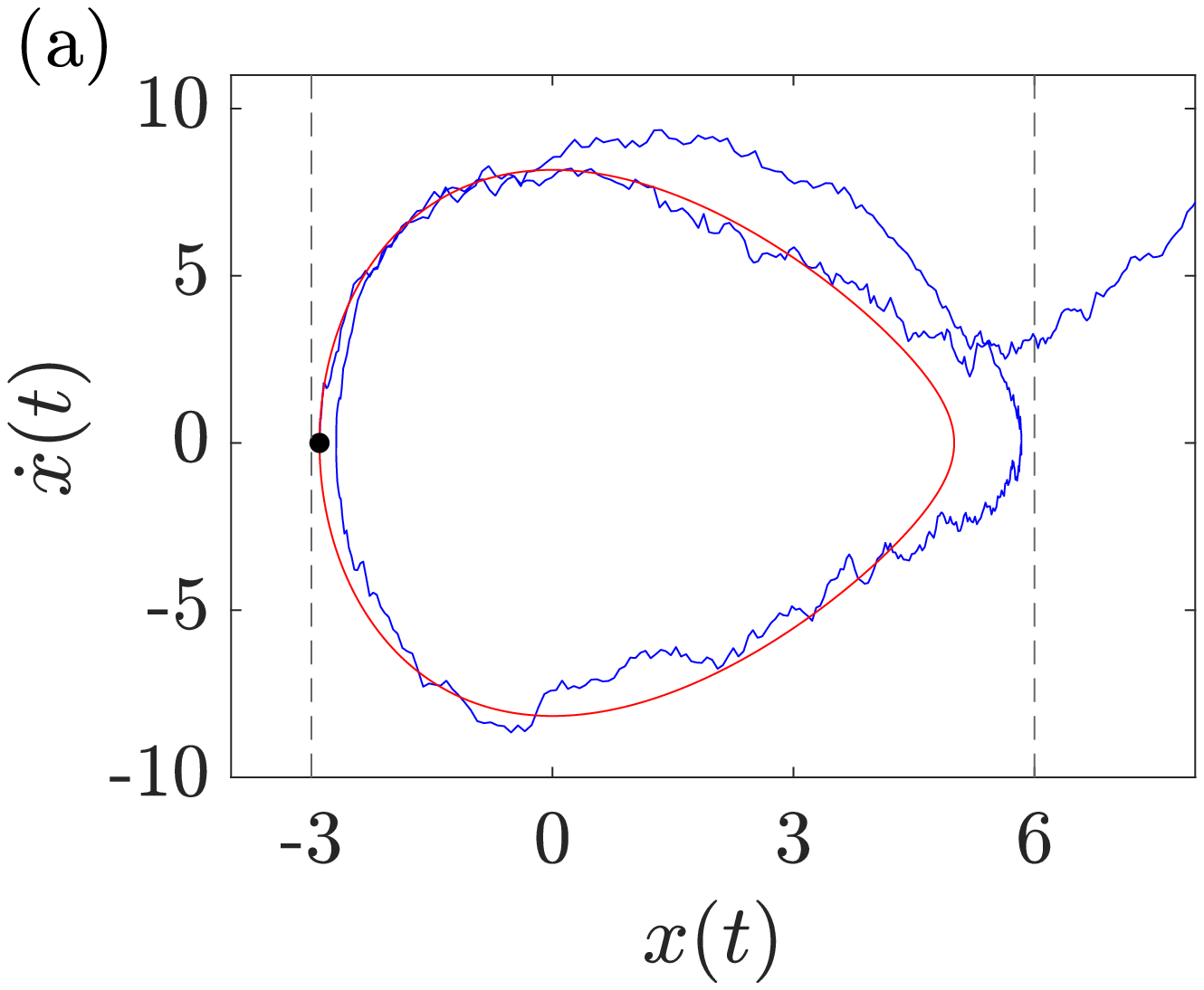}
		\includegraphics[width=0.45\textwidth]{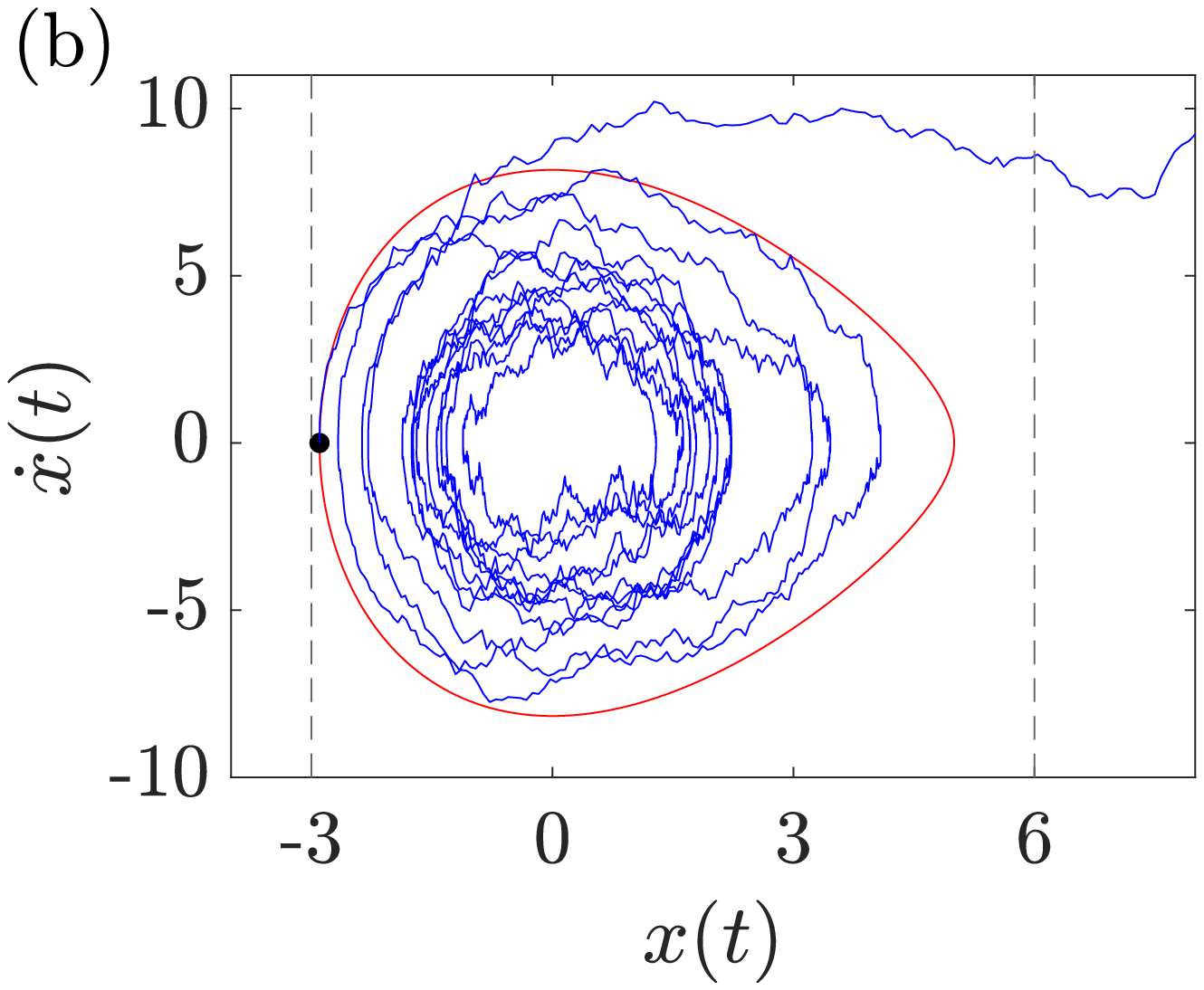}
	\end{center}
	\caption{Trajectories in phase space for two realizations of the process described by $ \ddot{x}+0.1 \dot{x}+6x-x^{2}= \xi(t) $. The black point marks the initial condition $ x_{0}=-2.899 $ and the limit cycle for the undamped noiseless system is drawn in red for reference. The vertical dotted lines denote the ends of the potential ($ x^{-} $ and $ x_{+} $). In (a), the trajectory completes one cycle before leaving the well, while in (b) the FPT is longer and it bounces back and forth multiple times before reaching $ x_{+} $.}
	\label{phase_space}
\end{figure}

Now, we present our numerical results for the Kramers problem without resetting. In the absence of noise and damping the system describes a periodical limit cycle $ L $, given that the initial condition is contained in the well, that is, $ -3<x_{0}<6 $. This attractor can be seen in Fig. \ref{phase_space} in red, the initial condition is chosen to be next to the left end, at $ x_{0}=-2.899 $, and it is marked with a black point. In the presence of damping, the origin is the only attractor and oscillations decay until this value. When noise is added, the particle follows a path similar to the limit cycle $ L $ in the sense that it bounces back and forth in the potential well, but periodicity is lost. Now the system can escape the potential well, once it crosses the line at $ x>6 $.

In Fig. \ref{phase_space} (a) the system completes a cycle before leaving the potential well. There are other trajectories for which a cycle is not completed and the particle never returns to the left side of the well, instead it immediately escapes. Finally, there are other trajectories for which this transient state lasts for longer periods of time and multiple cycles are described before leaving the well (see Fig. \ref{phase_space} (b)).

\subsection{MFPT distribution: noise and initial condition dependence}

The time that the particle takes on average to escape the potential barrier depends on the noise intensity and on the initial condition. Firstly, we computed the MFPT for a fixed initial condition and different noise levels to explore how much this parameter affects the escape times.

\begin{figure}[h]
	\begin{center}
		\includegraphics[width=0.45\textwidth]{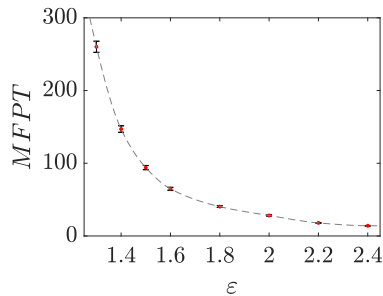}
	\end{center}
	\caption{MFPT dependence on the noise intensity for a fixed initial condition next to the left-hand side end of the well: $ x_{0}=-2.899 $. As expected, trajectories escape sooner for higher noise levels.}
	\label{MFPT_noise_dependence}
\end{figure}

In Fig. \ref{MFPT_noise_dependence}, we see the numerically calculated MFPT values for noise levels ranging from $ 3.6$ and $6.6 \% $ of the potential barrier and $x_{0}=-2.899 $. A strong noise, compared to the height of the potential barrier, makes trajectories escape faster than a small noise. Similar trends are found for other initial conditions. It can be seen that the error bars are smaller for stronger noises. This is because a strong noise also reduces substantially the standard deviation of the FPT distributions.

\begin{figure}[h]
	\begin{center}
		\includegraphics[width=0.45\textwidth]{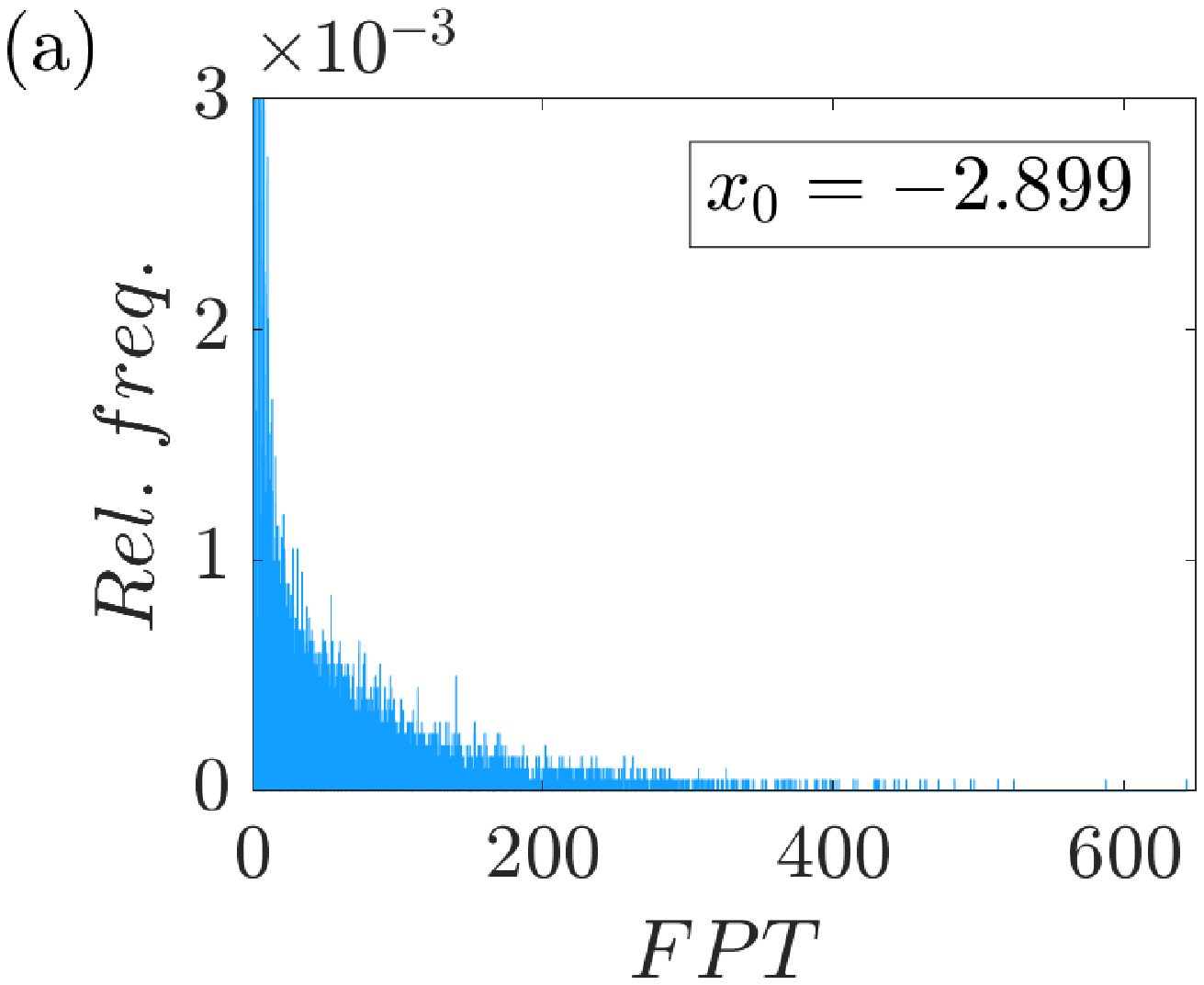}
		\includegraphics[width=0.45\textwidth]{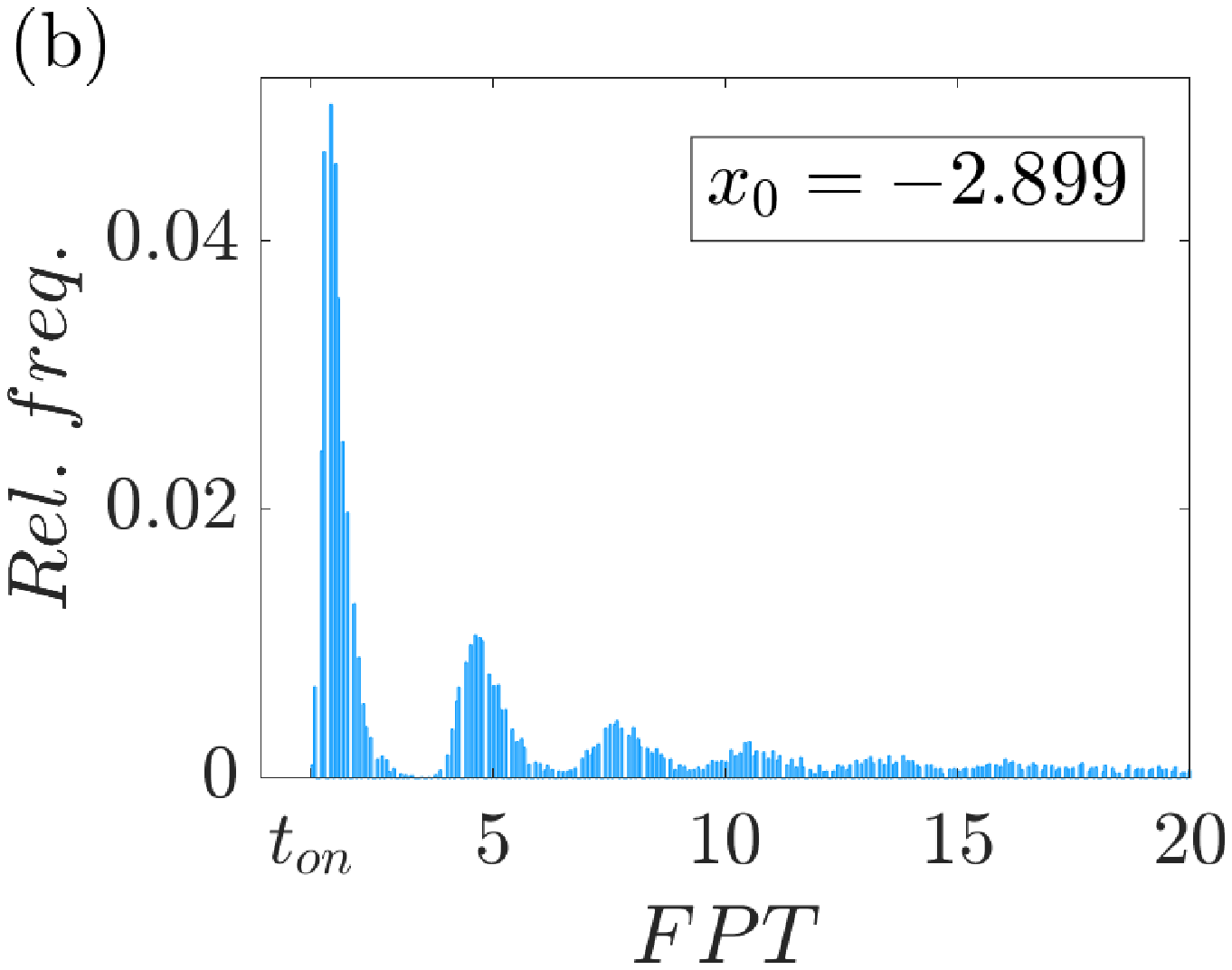}
		\includegraphics[width=0.45\textwidth]{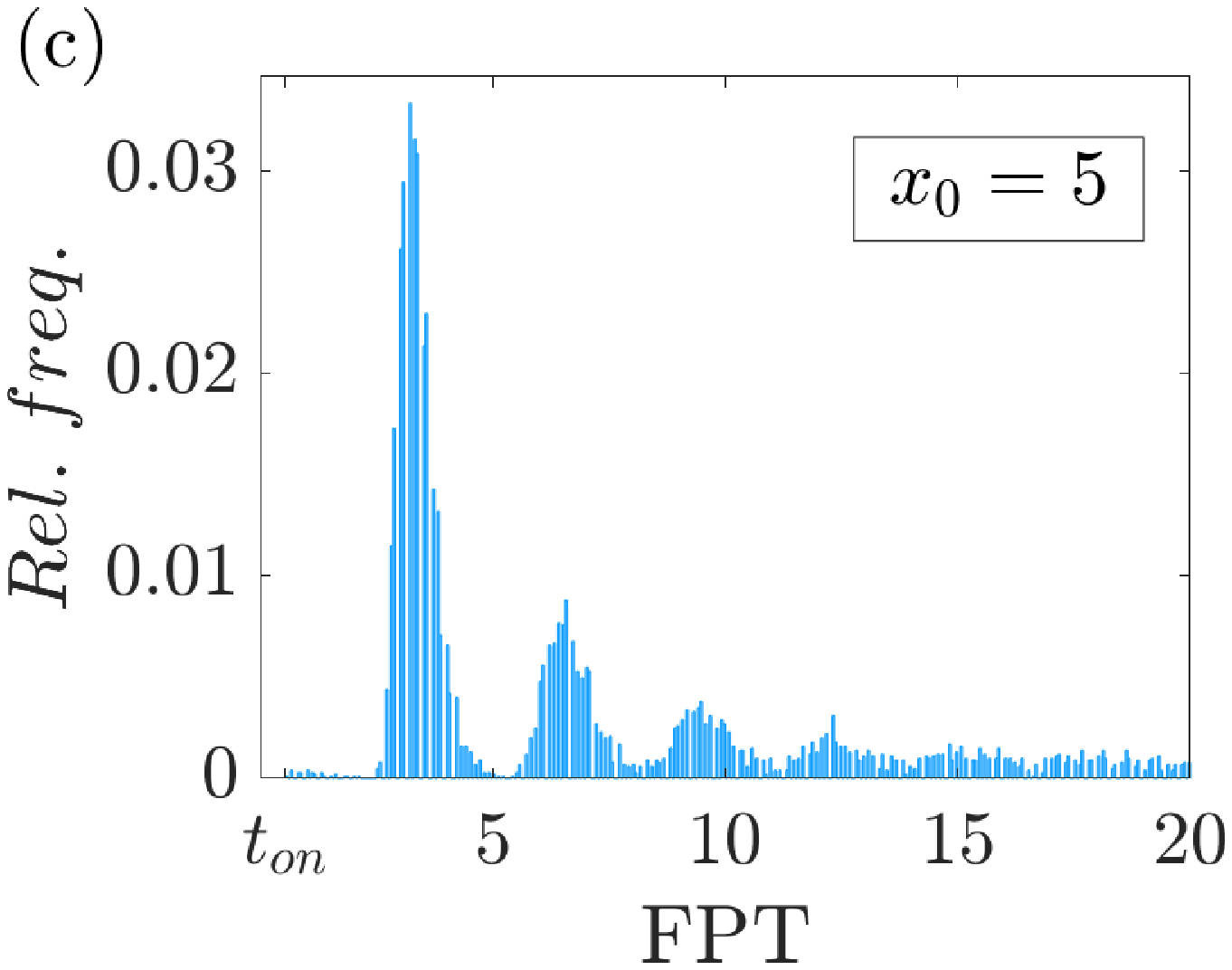}
		\includegraphics[width=0.45\textwidth]{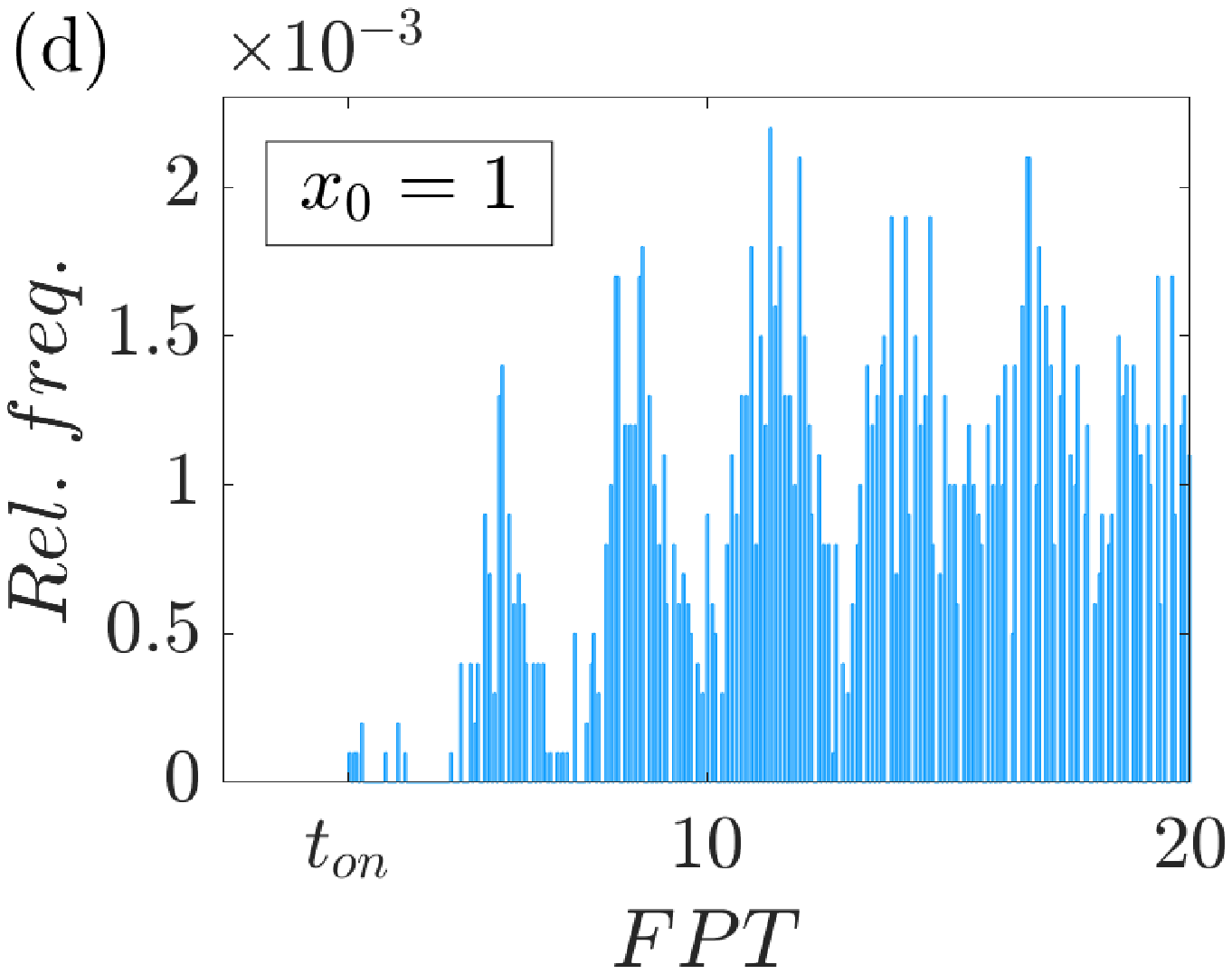}
	\end{center}
	\caption{Probability distributions for the first passage times. In (a) and (b) trajectories start at the left-hand side end of the well: $ x_{0}=-2.899 $. In (b) we show a zoom to visualize the wave pattern. In (c) trajectories start close to the escape: $ x_{0}=5 $ and the first wave, corresponding to escapes with bouncing, is very small. For trajectories far from the ends as $ x_{0}=1 $ in (d), the wave pattern is lost.}
	\label{histogram_complete_N1c8_x0_m2899}
\end{figure}

If we fix now the noise intensity to $ \varepsilon=1.8 $ ($ 5 \% $), we can explore the MFPT dependence with the initial condition. First of all, we fix $ x_{0} $ at one end of the potential well, $ x_{0}=-2.899 $. Starting on the left-hand side of the well, we show the first passage time distribution in Fig. \ref{histogram_complete_N1c8_x0_m2899} (a,b). In the $ y $-axis we represent the relative frequency, i.e., the empirical probability for a concrete value of FPT. The sum of the bar heights is equal to one.

In Fig. \ref{histogram_complete_N1c8_x0_m2899} (a), we see that the histogram is cut for a range of $ y $ values in order to visualize the bars for all FPTs. As it can be seen, most trajectories leave the well for small times. The last computed trajectory escapes at $ t=643.3 $. In Fig. \ref{histogram_complete_N1c8_x0_m2899} (b) we show the same histogram, now for FPTs ranging from $ 0 $ to $ 20 $. Approximately, $ 58 \% $ of the trajectories are in this range. From this perspective, we can appreciate a wave pattern. The reason is the following. From $ 0 $ to a time value called $ t_{on} $ no trajectories escape as $ t_{on} $ is the minimum time required for the particle to escape due to its finite velocity. Afterwards, there is a big wave corresponding to the trajectories that cross the well, passing through the origin and escape. The second wave corresponds to trajectories that bounce back and forth once as in Fig. \ref{phase_space} (a). The following waves correspond to an increasing number of swings. The waves rapidly become difficult to distinguish as the trajectories take different times to complete swings. Three long swings  (slow and/or big amplitude cycles) may take the same amount of time as four quick swings (fast and/or small amplitude cycles).

For other initial conditions, next to the ends of the potential well, similar FPT distributions are found. For instance, for $ x_{0}=5 $ (Fig. \ref{histogram_complete_N1c8_x0_m2899} (c)), the first wave is very small as the chances for the system to escape immediately and not bouncing are quite small. Due to the shape of the potential, the system tends to cross the origin and bounce around that value. In that case, the second wave is the biggest one and the followings decrease in magnitude. For more extreme initial condition values, the big waves (the first for left-hand positions and the second for right-hand positions) are even more important as it is more likely to escape sooner.

The picture changes for initial conditions far from the ends as in Fig. \ref{histogram_complete_N1c8_x0_m2899} (d). For them, the decaying amplitude wave pattern cannot be distinguished and trajectories take much more time to escape (compare the scale in the $ y- $axis). Fot $ t=20 $, only $ 15.04 \% $ of the trajectories have escaped the well. Remember that for $ x_{0}=-2.899 $, the $ 58 \% $ of the trajectories had escaped.

The shape of the FPT distributions without resetting plays a crucial role when resetting is added: it is the key to settle whether SR is beneficial or not. There is an interesting metric called coefficient of variation ($ CV $), which measures the ratio of the standard deviation and the mean ($ CV=\sigma / \tilde{\mu} $), that is, the relative fluctuation. If the values of the distribution are disperse with respect to the mean, then $ CV>1 $; if not, $ CV<1 $. Long tail distributions, like the ones found for initial conditions at the ends of the well, tend to have $ CV>1 $.

It has been shown that when $ CV>1 $, resetting helps the search process rendering smaller FPT values \cite{Pal2017}. If the fluctuations of the associated first-passage times about their mean are sufficiently strong, resetting tames some of the fluctuations. This can be easily understood looking at the histogram in Fig. \ref{histogram_complete_N1c8_x0_m2899} (a, b, c). Restarting the process would eliminate the long lasting trajectories that delay the overall process. After the system is completely reset, we find the same histogram once again. We let some fast trajectories escape and reset the process once again.

\begin{figure}[h]
	\begin{center}
		\includegraphics[width=0.45\textwidth ]{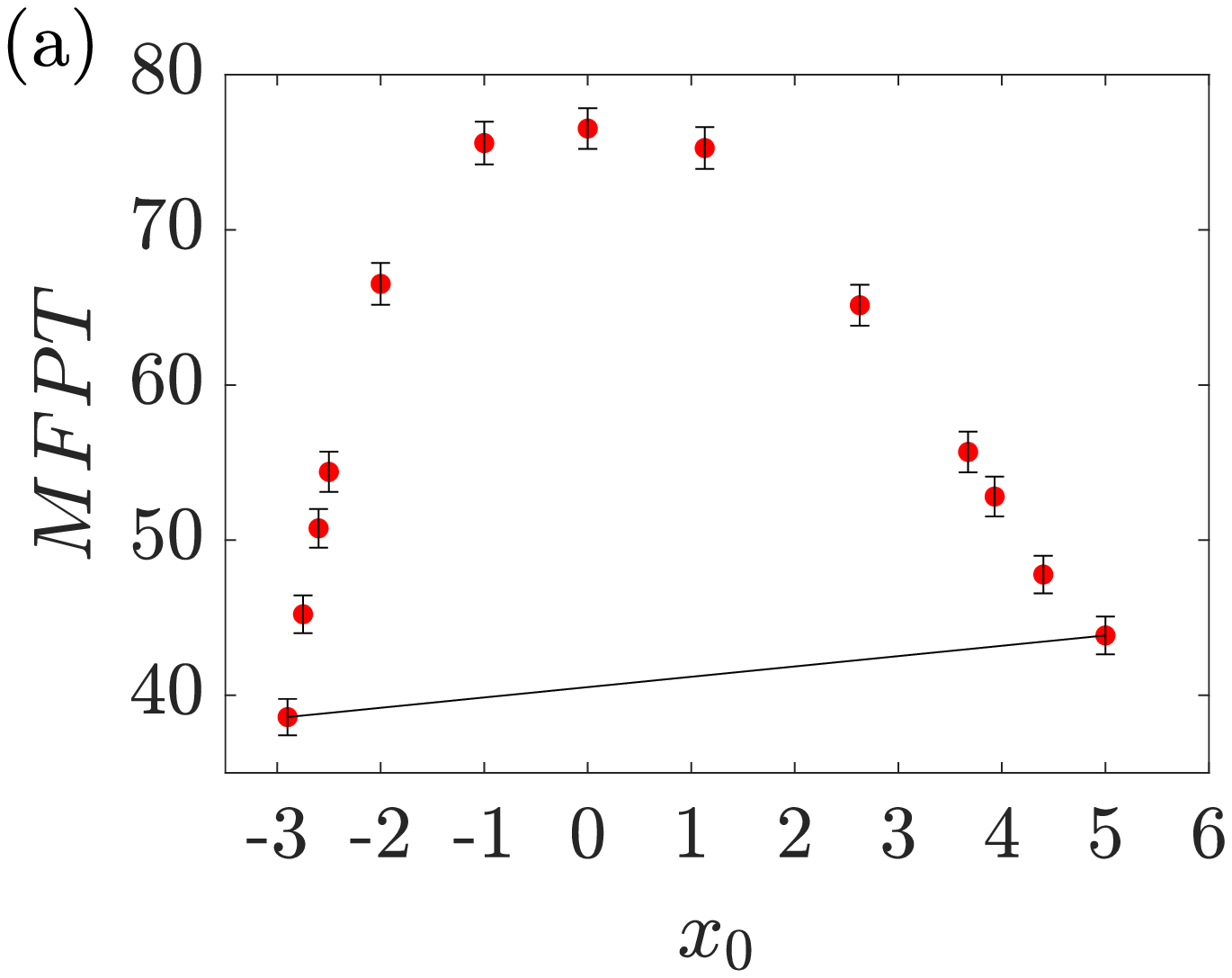}
		\includegraphics[width=0.45\textwidth]{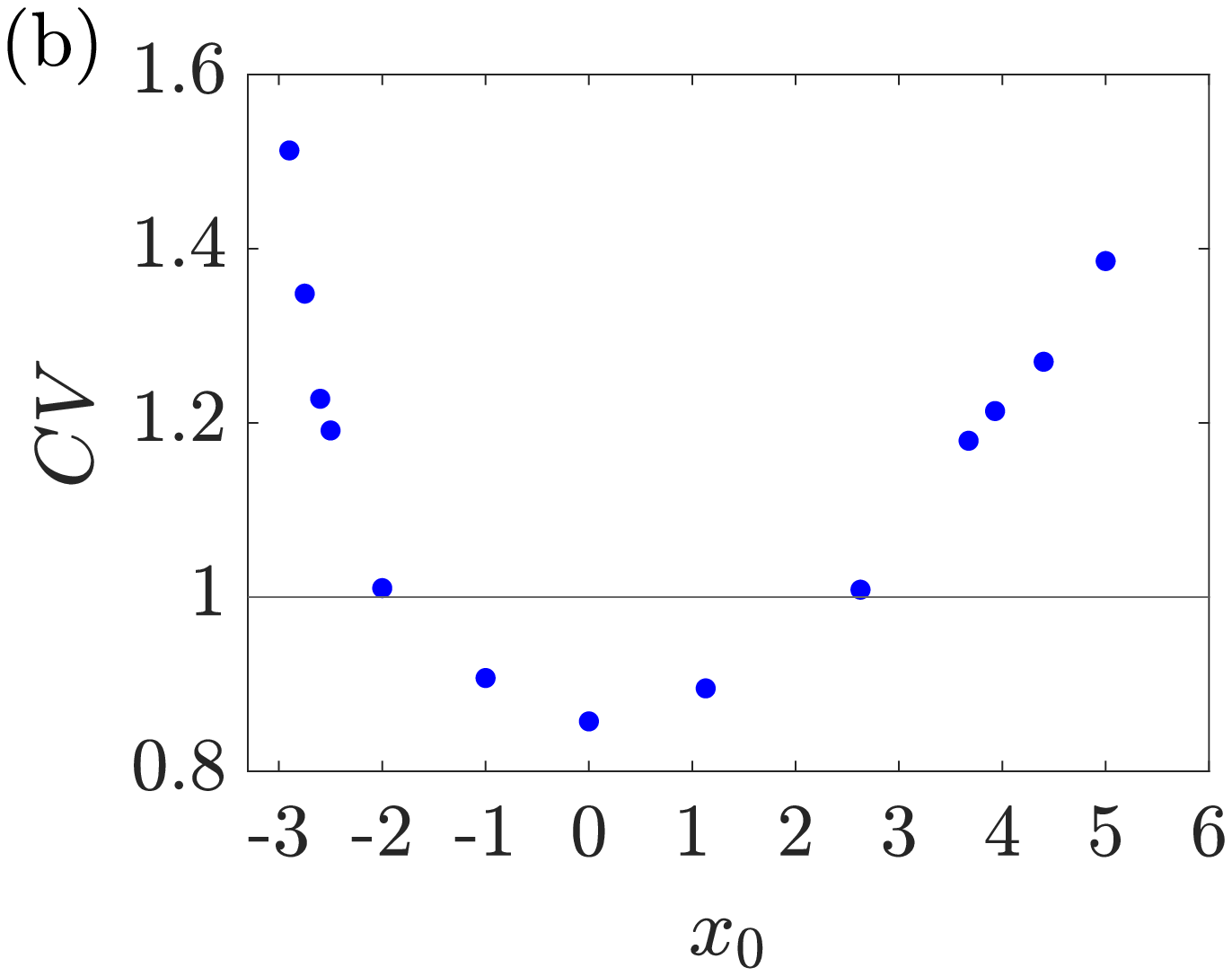}
	\end{center}
	\caption{Mean first passage time and coefficient of variation dependence on the initial condition. (a) Realizations starting close to the ends of the wells escape twice faster than those starting at the bottom of the potential. Due to the asymmetry of the potential, two initial conditions located at the same height (as the ones connected by a line in the figure), thus having the same initial energy, escape at different times. In (b) the $ CV $ serves as an indicator for resetting positions that would reduce the MFPT. The limit $CV=1 $ is achieved for positions at approximately a $ 40\% $ of $ E_{0} $.  }
	\label{MFPT_CV_IC}
\end{figure}

Figure \ref{MFPT_CV_IC} shows how the MFPT and $ CV $ depend on the initial condition. Points to the right and the left of the origin have the same initial energy. The MFPT is longer for trajectories starting around the origin, while trajectories next to the ends escape twice faster. It is important to notice that due to the asymmetrical shape of the potential, the MFPTs are not symmetrical at both sides of the origin. For visualization, we have included a line that connects two points for which the initial energy is the same. That is, $ \dot{x}^{2}/2=0 $ and $ V(x)=6 x_{0}^{2}/2- x_{0}^{3}/3 $ with $ x_{0}=(-2.899, 5) $. As it can be seen, points initially at the same energy show different MFPTs whether they are initially located further or closer to the escape at $ x_{+}=6 $. For instance, $ x_{0}=5 $, which is close to the escape, shows a longer MFPT because the first wave in the FPT distribution diagram is very small, as previously mentioned. This effect is more pronounced at the ends of the well. The differences between the MFPT at left- and right-handed locations diminish next to the origin.

In Fig. \ref{MFPT_CV_IC} (b), we see that the $ CV $ dependence on the initial condition is not symmetrical either. Furthermore, points above the line $ CV=1 $ mark initial conditions (and consequently reset positions as we are assuming it to be the same) for which resetting is beneficial, in the sense that it reduces the MFPT. These points are located at the ends of the well, while resetting to points near the origin would increase the MFPT. The limit is approximately at a energy level equal to $ 14.6 $ units of energy. Thus, resetting to points located (at least) at a $ 41 \% $ of the height of the potential barrier is needed for SR to be advantageous.

\section{Resetting at deterministic times} \label{Section_3}

For the rest of the paper, we focus on resetting to points for which $ CV>1 $ for a fixed noise intensity: $ \varepsilon=1.8 $. In this section, the resetting events take place at fixed time intervals, which is the simplest strategy to start with.

\begin{figure}[h]
	\begin{center}
		\includegraphics[width=0.45\textwidth]{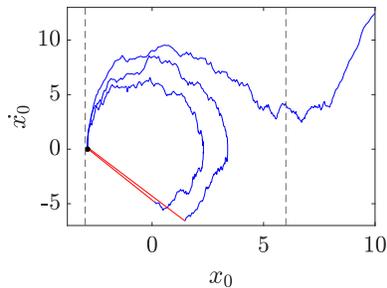}
	\end{center}
	\caption{Trajectory in phase space for a realization of the process described by $  \ddot{x}+0.1 \dot{x}+6x-x^{2}= \xi(t)  $, with resetting at fixed time intervals of $ t_{r}=2 $ (red). The black point marks the initial condition $ x_{0}=-2.899 $ and the vertical dotted lines denote the ends of the potential ($ x^{-} $ and $ x_{+} $). After two resetting events the particle leaves the well.}
	\label{resetting_phase_space}
\end{figure}

In Fig. \ref{resetting_phase_space} we show the phase space for a trajectory with resetting. The black point marks the initial condition, which is also the reset position. The trajectory moves around freely diffusing in the potential landscape until at $ t=2 $, the particle is immediately reset to the initial position. This bit is marked in red. Here, we are assuming that reset is instantaneous. Obviously, this is a simplification as in reality, resetting would always take some time. However, we assume that resetting is quick enough compared to the diffusion movement. For references in which stochastic resetting with time overheads is considered see \cite{Maso-Puigdellosas2019, Reuveni2016}; for refractory periods after the resetting, see \cite{Evans2019} and for stochastic returns using an external trap, see \cite{Gupta2020}. At $ t=4 $, the particle is reset once again. As it can be seen, after two resetting events, the system escapes the potential barrier.

Now, we calculate the FPT distribution for a constant time reset $ t_{r}=2 $. We simulate $ 10^{4} $ trajectories and measure the time they need to escape from the well. The numerical parameters are the same as in the previous section, except that we introduce now callback functions to account for the resetting events.

\begin{figure}[h]
	\begin{center}
		\includegraphics[width=0.45\textwidth ]{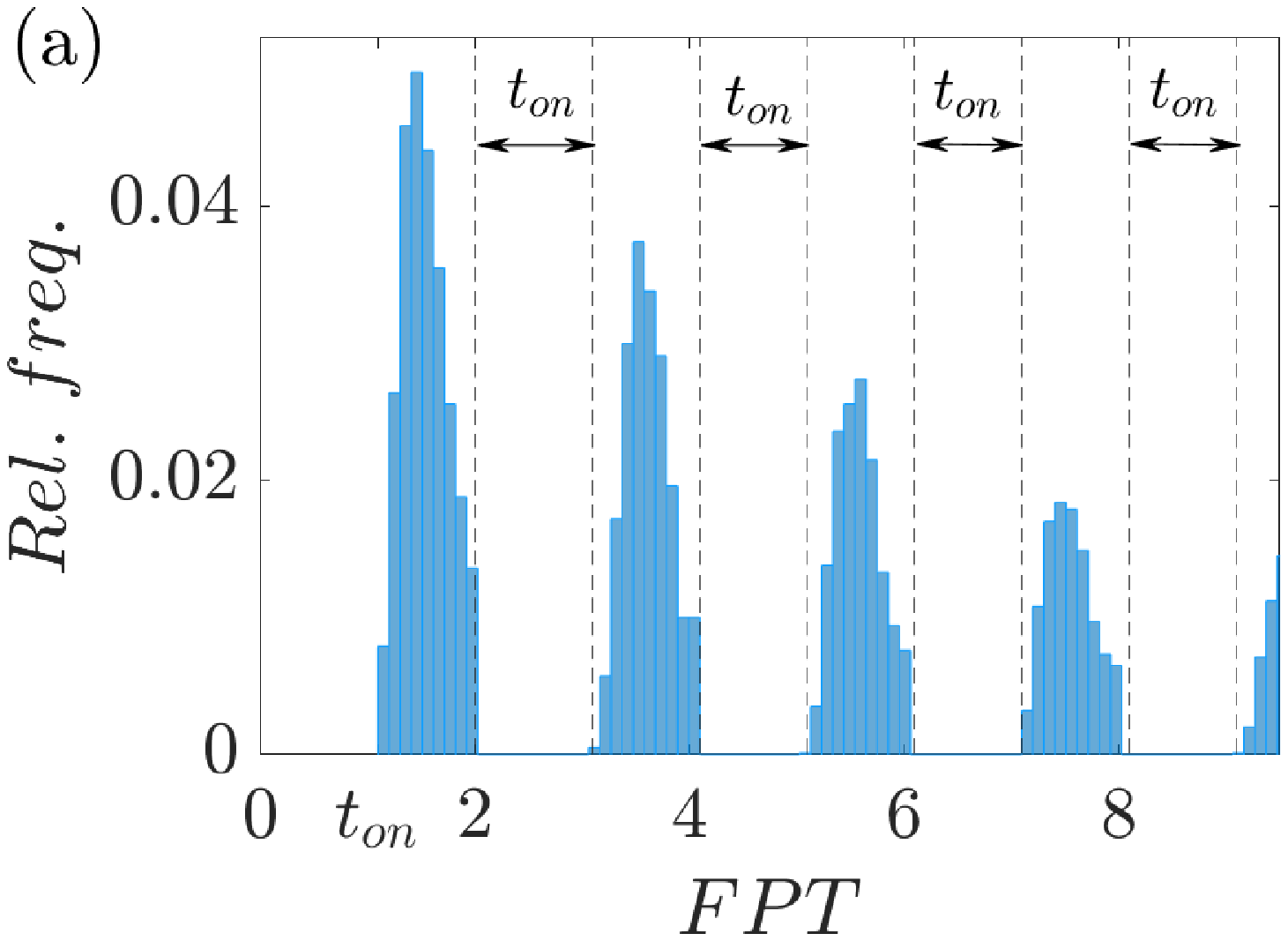}
		\includegraphics[width=0.45\textwidth]{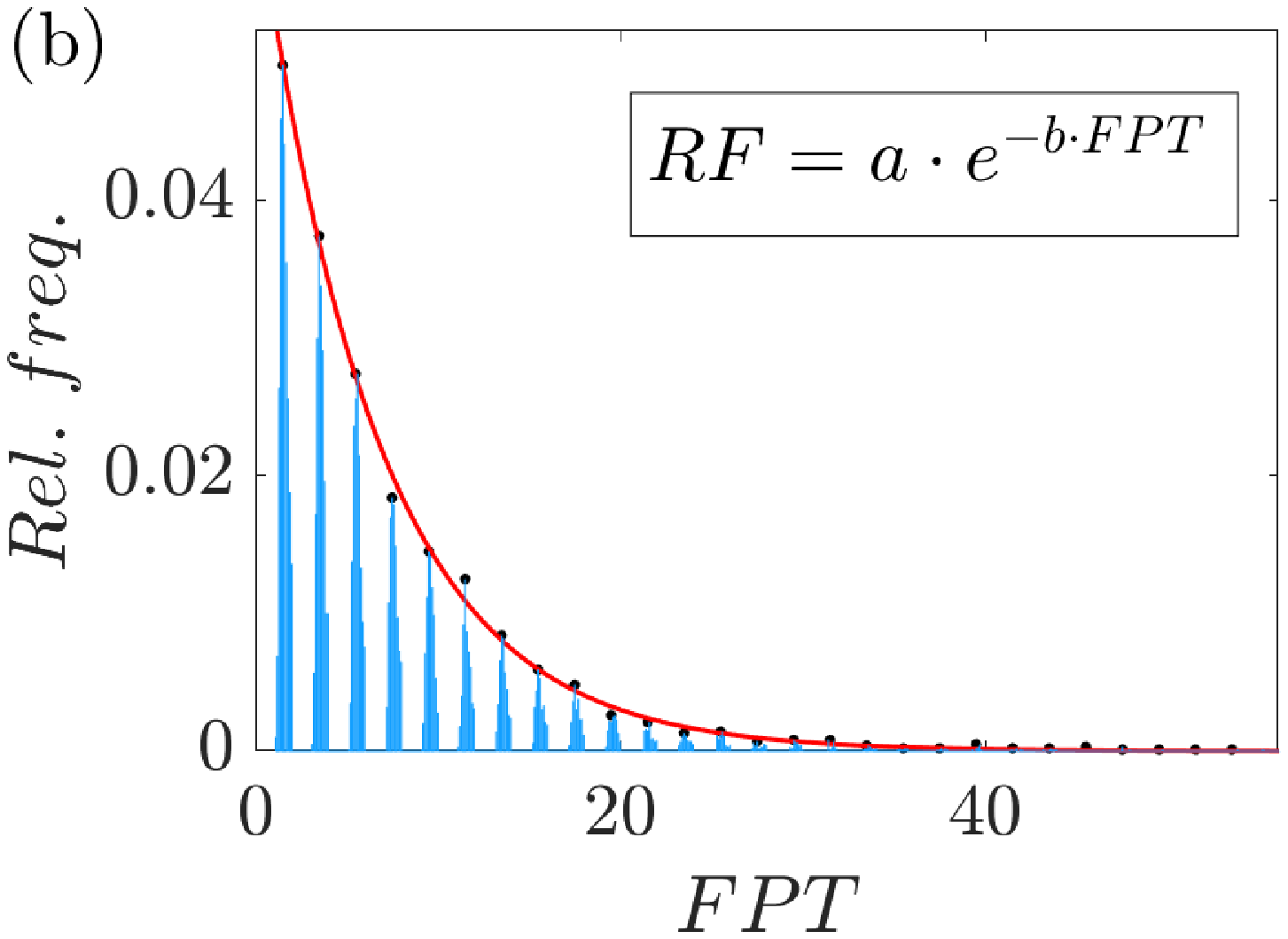}
	\end{center}
	\caption{Probability distribution for the first passage time with deterministic resetting for $ t_{r}=2 $. (a) After each resetting event no particles escape for a period of time equal to $ t_{on} $. (b) In the full view of the histogram it can be seen that the peaks of the waves decrease exponentially following the law shown in the inset.}
	\label{hist_SR}
\end{figure}

In Fig. \ref{hist_SR} (a) we see a zoom of the FPT distribution diagram to observe the wave pattern. From $ t=0 $ until the first resetting event, the diagram is the same as in the case without resetting (Fig. \ref{histogram_complete_N1c8_x0_m2899} (b)). From $ t=2 $ to $t=3.1  $ no trajectories escape the well as after the resetting, at least a time corresponding to $ t_{on}=1.1 $ has to pass. This fact produces a escape in the form of perfectly defined waves. In (b), we show the complete histogram where it can be seen that the peaks of each wave follow a exponential decay. The FPT of the trajectory that escapes the last is about $ 10 $ times smaller than without resetting. The relative frequency of the peaks ($ RF $) is adjusted to the following empirical trend:
\begin{equation}
	RF=a \cdot e^{-b \cdot FPT},
\end{equation}
where $ a=0.06274 $, $ b=0.1521 $ and $ R^{2}=0.9983 $. We conclude that stochastic resetting at fixed time intervals produces an escape in the form of waves of exponentially decreasing amplitude. This is logical as resetting happens for fixed time intervals and the number of remaining particles trapped in the well decreases at a fixed rate after each resetting event. Before the first reset, approximately the $ 27 \% $ of the trajectories have escaped. When the process is restarted we find the same escape pattern and again the $ 27 \% $ of the remaining trajectories escape. This leads to the exponential decay trend numerically observed in Fig. \ref{hist_SR} (b).

Now, we may ask ourselves if there is an optimum value for $ t_{r} $, this is, if resetting at a certain time interval is more beneficial than others. For that purpose, we calculated the MFPT for a range of $ t_{r} $ values, keeping the rest of parameters fixed.

\begin{figure}[h]
	\begin{center}
		\includegraphics[width=0.45\textwidth ]{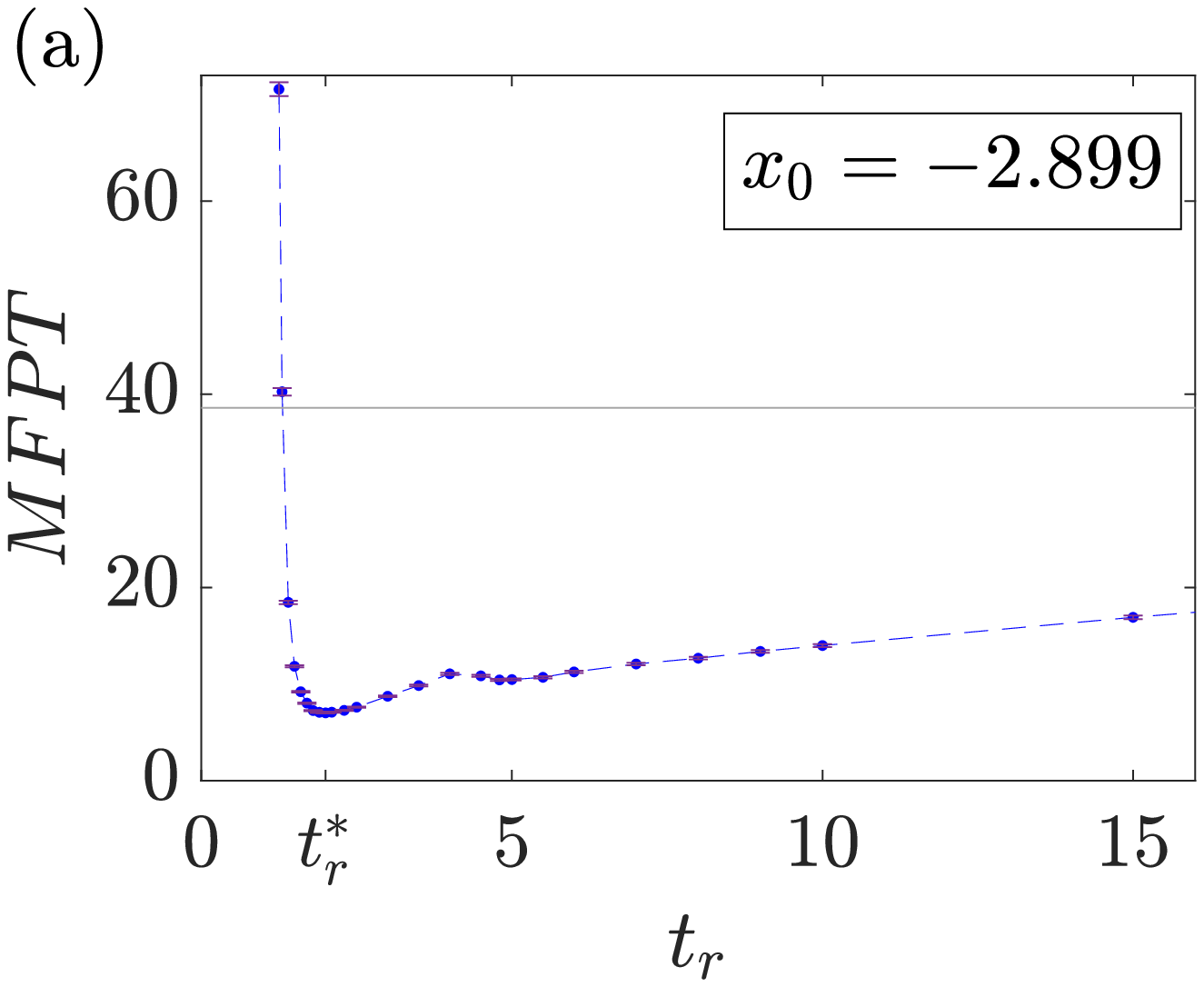}
		\includegraphics[width=0.45\textwidth]{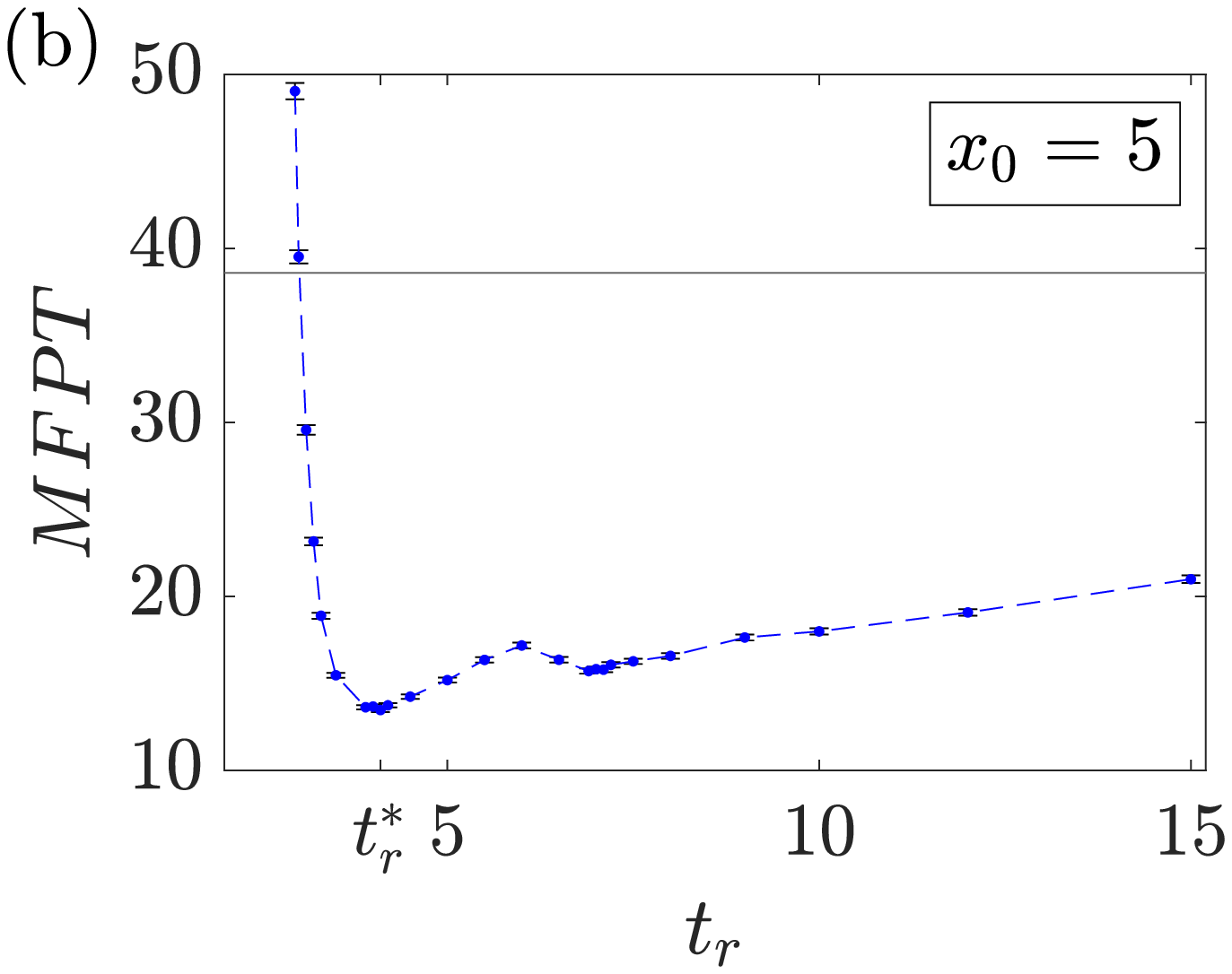}
	\end{center}
	\caption{Mean first passage time dependence on the resetting time interval $ t_{r} $ for an initial condition close to the left-hand end of the well (a) and close to the right-hand end of the well and the escape (b). In both cases there is an optimal value $ t_{r}^{*} $ for which the MFPT is minimum. This value corresponds to resetting after the big wave of trajectories in the non resetting case have escaped (see Fig. \ref{histogram_complete_N1c8_x0_m2899}). For shorter values of $ t_{r} $ SR may hinder the escape process. The horizontal line marks the MFPT for the case with no resetting for reference. Surprisingly, we find a second local minimum after $ t_{r}^{*} $, corresponding to the second biggest wave in the no resetting scenario. }
	\label{MFPT_tr}
\end{figure}

Figure \ref{MFPT_tr} shows this dependence for two different initial conditions located at the ends of the potential well ($ x_{0}=-2.899 $ and $ x_{0}=5 $). For these points the $ CV $ if higher, so resetting to them will be more effective. A horizontal line marks the MFPT without resetting for reference. It can be seen that in both cases for small values of $ t_{r} $, SR delays the escape process as the system is continuously returning to the initial point. Afterwards, we find a minimum which marks the optimal resetting time interval: $ t_{r}^{*} $. For $ t=t_{r}^{*} $, approximately a quarter of the trajectories have already left, so resetting gives a second chance to the trajectories that were stuck and did not escape fast enough. For higher values of $ t_{r} $, the MFPT tends slowly to the MFPT without resetting as resetting at very long intervals tends basically to the no resetting case. For other initial conditions with $ CV>1 $, but located further form the ends, the behavior is qualitatively the same but less drastic, that is, MFPT is not reduced as much.

Something interesting here is that SR is reduced for values close to the escape (on the right hand side of the well, figure (b)) and far from the escape (on the left hand side, figure (a)). This is in contrast with other papers as \cite{Singh2020} where they studied SR in a V-shaped potential and found that SR expedites the search process only when the reset location is on the same side as the absorbing trap. Here, the point $ x=6 $, that is, the top of the potential, is found more rapidly even if the reset location is on the contrary side. In fact, for $ x_{0}=-2.899 $, the MFPT reduces at its best to a $ 18\% $ of the no resetting MFPT value, while for $ x_{0}=5 $, the MFPT reduces to a $ 35\%  $.

Another difference that we find for our potential is the presence of a second local minimum in the MFPT. For $ x_{0}=-2.899 $, the global minimum is $ t_{r}^{*}=2 $ and the second minimum is found for $ t_{r}=4.8 $. The explanation for this is the following. The first minimum corresponds to resetting at a time for which in the no resetting case the peak of the first wave has just passed. This means that resetting right after the first bunch of trajectories has escaped is the best strategy. The second minimum corresponds to resetting right after the peak of the second wave in the no resetting case has passed. This is more beneficial than other values right before or afterwards, but still the best strategy is resetting after the first wave. Other minima are expected after the following waves but they cannot be appreciated at our resolution.

In the case of $ x_{0}=5 $, $ t_{r}^{*}$ is found for a later value than in the previous case. Concretely, $ t_{r}^{*}=4.1$ because in the no resetting situation the important wave was the second one and the first one could almost not be noticed. The second minimum is found at $ t_{r}=7.1 $ and corresponds to third wave in the case without resetting.

These results show that in a Kramers problem the best resetting strategy can be designed just by studying the case without resetting and fixing resetting times after the first big escape wave.

\section{Resetting at stochastic times} \label{Section_4}

The nature of the resetting can be diverse. In the literature of SR it is mostly considered to be a Poisson event. This is because they are ubiquitous in real world problems and occur independently at a constant average rate, $ r $. They are memoryless processes in the sence that after an event has taken place, the probability of another event taking place afterwards is not increased nor decreased. Other resetting time configurations considered in preceding works are deterministic \cite{Bhat2016}, power-law distributed \cite{Nagar2016} or show a time \cite{Pal2016} or position dependence \cite{Evans2011a}.

\begin{figure}[h]
	\begin{center}
		\includegraphics[width=0.45\textwidth]{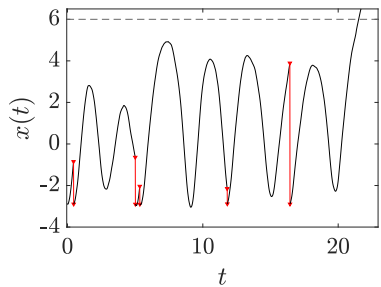}
	\end{center}
	\caption{Time series for a realization of the process defined by $ \ddot{x}+0.1 \dot{x}+6x-x^{2}= \xi(t) $, with resetting at exponentially distributed times. After five resetting events (red), the trajectory crosses the horizontal dotted line that marks the top of the potential and escapes. }
	\label{timeseries_poisson}
\end{figure}

For Poisson events, the time that passes between successive events follows an exponential distribution. Its probability density function is of the form $ f(t;r)=r \cdot e^{-r \cdot t} $, where $ r $ is the rate parameter and the mean of the distribution is $ 1/r $. In Fig. \ref{timeseries_poisson}, we show a time series of the process with SR for a given value of $ r $. The initial condition and resetting point is $ x_{0}=-2.899 $. Between the resetting events, marked in red, the system is diffusing in the potential well following Eq. \ref{Langevin}. For this realization, five resetting events occur before the particle leaves the well trespassing $ x=6 $.

\begin{figure}[h]
	\begin{center}
		\includegraphics[width=0.45\textwidth]{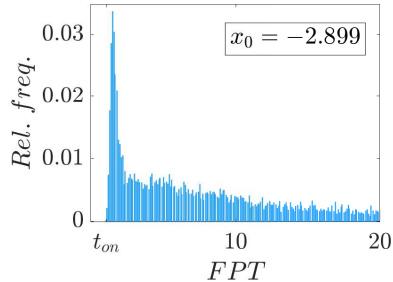}
	\end{center}
	\caption{Probability distribution for the first passage times with stochastic resetting for $ x_{0}=-2.899 $. In this case, the wave pattern is lost as each realization has different resetting times. For other resetting positions for which $ CV>1 $ we find qualitatively similar diagrams.}
	\label{histogram_Poisson}
\end{figure}

For a large number of realizations, we get a histogram like the one in Fig. \ref{histogram_Poisson}. In this case, the wave pattern found for deterministic resetting is lost as each realization has different resetting times. Comparing this with Fig. \ref{histogram_complete_N1c8_x0_m2899} (b), we can see that the height of the main peak is reduced as some trajectories return to the initial point before the optimal value, delaying the escape process. However, the MFPT and the standard deviation are largely reduced. Until $ t=20 $, approximately the $ 85 \% $ of the trajectories have escaped, while for the case without resetting it was only the $ 58 \% $.

\begin{figure}[h]
	\begin{center}
		\includegraphics[width=0.45\textwidth]{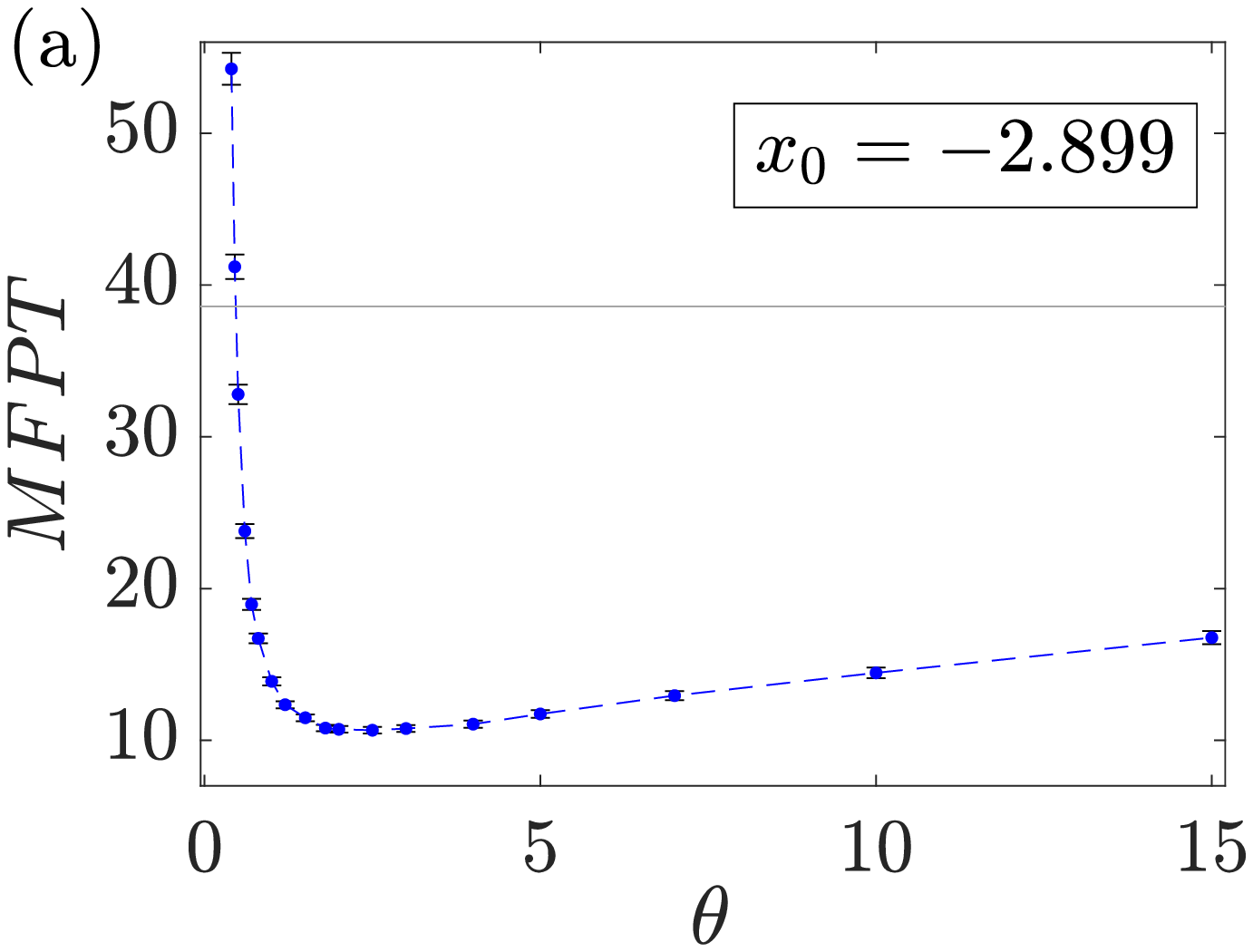}
		\includegraphics[width=0.45\textwidth]{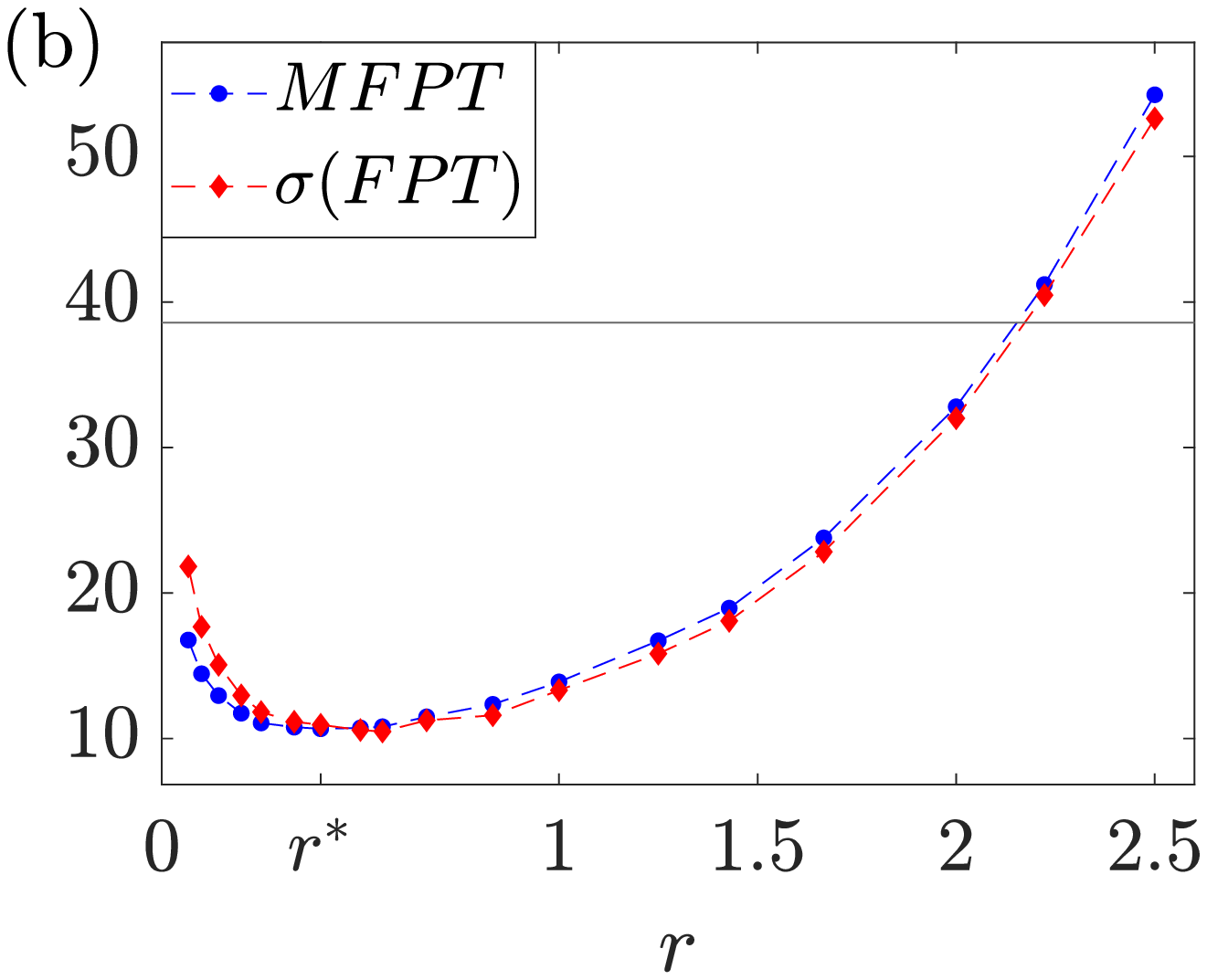}
	\end{center}
	\caption{Mean first passage time dependence on the resetting rate. (a) The MFPT is represented against $ \theta=1/r $. The minimum corresponds to a mean resetting interval slightly above the deterministic resetting $ t_{r}^{*} $ value. (b) The MFPT and the standard deviation are equal for the optimal rate, $ r=r^{*}=0.4 $.}
	\label{MFPT_poisson}
\end{figure}

As in the previous section, we are interested in finding the optimal resetting rate for which the MFPT is minimum. For a given initial condition $ x_{0}=-2.899 $, we numerically calculate the escape time for $ 10^{4} $ realizations of the process. In Fig. \ref{MFPT_poisson} (a), we show how the MFPT varies with $ \theta=1/r $. The MFPT without resetting is marked as a horizontal line for reference. We find a similar behavior as in the deterministic resetting scheme. Only that now the second local minimum is not present. This is due to the fact that resetting times are not always the same and, consequently, they do not correspond exactly with any of the waves of the no resetting case.

The optimal rate is found for $ \theta=2.5 $, which corresponds to a mean resetting interval of $ 2.5 $. This is slightly above the deterministic resetting $ t_{r}^{*} $ value. For practical purposes any resetting rate that arises a mean resetting interval similar to the time for the first wave in the case without resetting, will be beneficial. In the case of resetting to points on the right-hand side of the well, then we are referring to the second wave which is the big one. For $ x_{0}=5 $, the optimal resetting interval is $ \theta=5.5$, which is also slightly above the deterministic resetting case.

Finally, if the reset position is far from the ends, this is, points for which $ CV \approx 1 $, then the optimal $ \theta $ increases and no longer corresponds to any of the waves. This is because the first waves have a small amplitude. Consequently, the MFPT is not reduced as much. For instance, if $ x_{0}=-2 $ for which $ CV=1.0103 $, then resetting is found for $ \theta \sim 35 $ and the MFPT reduces only to a $ 91 \% $ of the MFPT without resetting. For these reset positions, the standard deviation barely decreases for $ r^{*} $ as opposed to reset positions at the ends of the wells (see Fig. \ref{MFPT_poisson} (b)). In that case, the statistical fluctuations due to the Monte Carlo method are a limitation as they depend on the standard deviation, remember Eq. \ref{MC_error}. To properly study these reset positions far from the ends, one would need to increase $ N $ by two or three orders of magnitude.

In Fig. \ref{MFPT_poisson} (b), we show the MFPT (blue circles) and the standard deviation (red diamonds) dependence with $ r $. Here, we obtain an analytical result derived by Reuveny \cite{Reuveni2016}, which asserts that the relative fluctuation ($ CV $) in the FPT of an optimally restarted process is always unity. That is, for $ r^{*} $:
\begin{equation}
	MFPT =\sigma(FPT).
\end{equation}

This result holds for any FPT process, given that the restart time is taken from a exponential distribution (constant restart rate).

At $ r^{*} $, the MFPT reduces approximately to a $ 28 \% $ of the MFPT without resetting. Comparing this to the deterministic resetting strategy, the Poisson resetting is less effective.

\section{Conclusions} \label{Section_5}

In this paper we have analyzed the stochastic resetting strategy in a Kramers problem. We applied this technique to a particle diffusing in a potential with one escape in order to minimize the mean escape time. We have shown that stochastic resetting may be beneficial even for restarting positions on the contrary side of the escape.

Following an analytical approach is not always possible as in the case of complex potentials because no closed-forms expressions arise. We have demonstrated that the Monte Carlo approach, which is needed for realistic situations, is in agreement with some analytical results and that this technique allows to explore the system in a different way.

We have found that in the case without resetting the escape times distribution forms a wave pattern. As expected, the MFPT was found to decrease for increasing noise levels. Also, it was found to decrease for trajectories starting at a higher energy level (next to the ends of the well). Surprisingly, we compared two trajectories starting with the same energy level and we observed that the ones further from the escape presented lower escape times. The numerically calculated coefficient of variation also served as an indicator for the efficiency of the stochastic resetting strategy. We found that $ CV>1 $ for trajectories starting at least at a $ 41 \% $ of the height of the potential barrier.

Then, we applied two resetting strategies: resetting at deterministic times and resetting at exponentially distributed times (treating the resetting events as Poisson events). In the first case, resetting drastically reduces the MFPT and the escape times distribution takes the form of waves with exponentially decreasing height. We also found that the best resetting strategy can be predicted in advance as the optimal resetting time interval is related to the wave pattern in the no resetting case.

Resetting at exponentially distributed times is also beneficial, but is less efficient than a periodic resetting. The optimum was found for rates slightly smaller, that is, slightly longer mean time intervals than in the deterministic resetting situation. In this case we retrieve another important analytical prediction that asserts that for the optimal resetting rate, the MFPT equals the standard deviation. In both cases, resetting to the contrary side of the escape was the optimum.

Finally, we hope that our results and approach can be applied to areas where diffusion in a potential barrier is relevant such as atomic diffusion, chemical reactions or Josephson junctions, where SR can be an induced or natural process.

\bibliographystyle{elsarticle-num}
\bibliography{library}
\end{document}